\documentclass[reqno,12pt]{amsart}
\usepackage{amsthm}
\usepackage{amsmath,mathrsfs}
\usepackage{color}
\usepackage{amssymb}
\usepackage{hyperref}
\usepackage{enumerate}
\usepackage{latexsym,array}
\usepackage{amsfonts}
\usepackage{shadow}

\newtheorem{Pa}{Paper}[section]
\newtheorem{Tm}[Pa]{{\bf Theorem}}
\newtheorem{La}[Pa]{{\bf Lemma}}
\newtheorem{Dn}[Pa]{{\bf Definition}}

\newtheorem{Rk}[Pa]{{\bf Remark}}

\def\C{\mathbb C}
\def\hh{\mathbb{H}}
\author[D. Alpay]{Daniel Alpay}
\address{(DA) Department of Mathematics\\
Ben-Gurion University of the Negev\\
Beer-Sheva 84105 Israel} \email{dany@math.bgu.ac.il}

\author[F. Colombo]{Fabrizio Colombo}
\address{(FC) Politecnico di
Milano\\Dipartimento di Matematica\\Via E. Bonardi, 9\\20133
Milano, Italy}
\email{fabrizio.colombo@polimi.it}

\author[D. P. Kimsey]{David P. Kimsey}
\address{(DPK)Department of Mathematics\\
Ben-Gurion University of the Negev\\
Beer-Sheva 84105 Israel}
\email{dpkimsey@gmail.com}

\author[I. Sabadini]{Irene Sabadini}
\address{(IS) Politecnico di
Milano\\Dipartimento di Matematica\\Via E. Bonardi, 9\\20133
Milano, Italy}
\email{irene.sabadini@polimi.it}

\title[The spectral theorem for unitary operators based on the $S$-spectrum]
{The spectral theorem for unitary operators\\ based on the $S$-spectrum} \oddsidemargin
0.2in \evensidemargin 0.2in \topmargin -0.5in \textwidth
15.5truecm \textheight 23truecm
\def\H{\mathbb H}

\def\R{\mathbb R}

\def\C{\mathbb C}

\def\(s){\mathscr S(\R\times\R)}

\newcommand{\ZZ}{\mathbb{Z}}
\newcommand{\CC}{\mathbb{C}}
\newcommand{\HH}{\mathbb{H}}
\newcommand{\cH}{\mathcal{H}}
\newcommand{\cB}{\mathbf{B}}



\thanks{D. Alpay thanks the Earl Katz family for endowing the chair
which supported his research. D. P. Kimsey gratefully acknowledges the support of a Kreitman postdoctoral fellowship. F. Colombo and I. Sabadini
acknowledge the Center for Advanced Studies of the Mathematical
Department of the Ben-Gurion University of the Negev for the
support and the kind hospitality during the period in which part
of this paper has been written.}

\begin{document}
\begin{abstract}
The quaternionic spectral theorem has already been considered in the literature, see e.g. \cite{14}, \cite{sc}, \cite{Viswanath}, however, except for the finite dimensional case in which the notion of spectrum is associated to an eigenvalue problem, see \cite{fp}, it is not specified which notion of spectrum underlies the theorem.

 In this paper we prove the quaternionic spectral theorem for unitary operators using the $S$-spectrum.  In the case of quaternionic matrices, the $S$-spectrum coincides with the right-spectrum and so our result recovers the well known theorem for matrices.
  The notion of $S$-spectrum is relatively new, see \cite{book_functional}, and has been used for quaternionic linear operators, as well as for $n$-tuples of not necessarily commuting operators, to define and study a noncommutative versions of the Riesz-Dunford functional calculus.\\
 The main tools to prove the spectral theorem for unitary operators are the quaternionic version of Herglotz's theorem, which relies on the new notion of $q$-positive measure, and quaternionic spectral measures, which are related to the quaternionic Riesz projectors defined by means of the $S$-resolvent operator and the $S$-spectrum.
\\
 The results in this paper restore the analogy with the complex case in which the classical notion of spectrum appears in the Riesz-Dunford functional calculus as well as in the spectral theorem.
\end{abstract}

\maketitle
\tableofcontents
\parindent 0cm

\section{Introduction}
\setcounter{equation}{0}

One of the main motivations to study spectral theory of linear operators in the quaternionic setting is due to the fact that Birkhoff and von Neumann, see \cite{BvN}, showed that there are essentially two possible settings in which to write the Schr\"odinger equation: one may use complex-valued functions or one may use quaternion-valued functions. Since then, many efforts have been made by several authors, see \cite{adler, 12, 14, 21}, to develop a quaternionic version of quantum mechanics.

Fundamental tools in this framework are the theory of quaternionic groups and semigroups on quaternionic Banach spaces which have been studied only recently in the papers
\cite{perturbation, evolution, GR} using the notion of $S$-spectrum and of $S$-resolvent operator as well as the spectral theorem, which is the main result of this paper.

To fully understand the aim of this work, we start by recalling some basic facts in complex spectral theory.
Let $A$ be a linear  operator acting on a complex Banach space $X$, and let $\sigma(A)$ and $\rho(A)$  be the spectrum and the resolvent sets of $A$, respectively.
One of the most natural ways to associate to a linear operator $A$ the linear operator $f(A)$ is to use the Cauchy formula for holomorphic functions
$$
f(A)=\frac{1}{2\pi i}\int_{\partial\Omega} (\lambda I-A)^{-1} f(\lambda) d\lambda,
$$
where $\partial\Omega$ is a smooth closed curve that belongs to the resolvent set of $A$ and $f$ a holomorphic function on an open set $\Omega$ which contains the spectrum of $A$.
This holomorphic functional calculus is known as Riesz-Dunford functional calculus, see \cite{ds}.\\
To any linear operator $A$, it is possible to associate the notion of spectral measures, which can be written explicitly using the Riesz-projectors, as described below.
A subset of $\sigma(A)$ that is open and closed in the relative topology of $\sigma(A)$ is called a spectral set.
The spectral sets form a Boolean algebra and with each spectral set $\sigma$ one can associate the projector operator
$$
P(\sigma)=\frac{1}{2\pi i}\int_{C_\sigma} (\lambda I-A)^{-1}d\lambda
$$
where $C_\sigma$ is a smooth closed curve belonging to the resolvent set $\rho(A)$, such that $C_\sigma$ surrounds $\sigma$ but no other points of the spectrum.
A spectral measure in the complex Banach space $X$ is then a homomorphic map of the Boolean algebra of the sets into the Boolean algebra of projection operators
in $X$ which has the additional property that it maps the unit in its domain into the identity operator in its range.

As is well known, the spectrum $\sigma(A)$ appearing in the definition of the  Riesz-projectors $P(\sigma)$
is the same spectrum on which is supported the spectral measure $E(\lambda)$  appearing in
 the spectral theorem for normal linear operators in a complex Hilbert space.
Precisely, for a normal  linear operator $B$ on a Hilbert space, given a continuous function $g$ on the spectrum $\sigma(B)$, we have
$$
g(B)=\int_{\sigma(B)}g(\lambda) dE(\lambda).
$$
In the quaternionic setting, before the introduction of the $S$-spectrum, see \cite{book_functional}, two spectral problems were considered. We discuss the case of
a right linear quaternionic operator (the case of a left linear operator being similar)  $T:\mathcal{V}\to \mathcal{V}$ acting on a quaternionic  two sided Banach space $\mathcal{V}$, that is $T(w_1\alpha+w_2\beta)=T(w_1)\alpha +T(w_2)\beta$, for
$\alpha,\beta\in\mathbb{H}$ and  $w_1$, $w_2\in \mathcal{V}$. The symbol $\mathcal{B}^R(\mathcal{V})$ denotes the Banach space of bounded right linear operators.
\\
The left spectrum $\sigma_L(T)$ of $T$ is related to the resolvent operator
$(s\mathcal{I}-T)^{-1}$
that is
$$
\sigma_L(T)=\{ s\in \mathbb{H}\ \ :\ \ s\mathcal{I}-T\ \ \ {\rm is\ not\  invertible\ in\ }\mathcal{B}^R(\mathcal{V}) \},
$$
where the notation $s\mathcal{I}$  in $\mathcal{B}^R(\mathcal{V})$ means that
$(s\mathcal{I} )(v)=sv$.
\\ The right spectrum $\sigma_R(T)$ of $T$ is associated with the  right eigenvalue problem, i.e. the search for nonzero vectors
satisfying $T(v)=vs$. It is important to note that
if $s$ is an eigenvalue, then all quaternions belonging to the sphere
$r^{-1}s r$, $r\in\mathbb{H}\setminus \{0\}$, are also
eigenvalues.
But observe that the operator  $ \mathcal{I}s -T$  associated to the right eigenvalue problem is not linear, so it is not clear what is the resolvent operator to be considered.

 A natural notion of spectrum that arises in the definition of the  quaternionic functional calculus is the one of $S$-spectrum.
In the case of matrices, the $S$-spectrum coincides with the set of right eigenvalues; in the general case of a linear operator, the point $S$-spectrum coincides with the set of right eigenvalues.

In the literature there are several papers on the quaternionic spectral theorem, see e.g. \cite{14,Viswanath}, however the notion of spectrum in use is not made clear.
Recently, there has been a resurgence of interest in this topic, see \cite{spectcomp}, where the authors prove the spectral theorem, based on the $S$-spectrum,
 for compact normal operators on a quaternionic Hilbert space.
In this paper we prove the quaternionic spectral theorem for unitary operators using the $S$-spectrum, which is then realized to be the correct notion of spectrum for the quaternionic spectral theory of unitary operators.

The  $S$-spectrum, see \cite{book_functional}, is defined as
$$
\sigma_S(T)=\{ s\in \mathbb{H}\ \ :\ \ T^2-2 {\rm Re}(s)T+|s|^2\mathcal{I}\ \ \
{\rm is\ not\  invertible}\},
$$
while the  $S$-resolvent set is
$$
\rho_S(T):= \mathbb{H}\setminus\sigma_S(T)
$$
where $s=s_0+s_1i+s_2j+s_3k$ is a quaternion,  $i$, $j$ and $k$ are the imaginary units of the quaternion $s$,  ${\rm Re}(s)=s_0$ is the real part and the norm $|s|$ is such that $|s|^2=s_0^2+s_1^2+s_2^2+s_3^2$.
Due to the noncommutativity of the quaternions, there are two resolvent operators associated with a quaternionic linear operator:
the left  and the right $S$-resolvent operators which are defined as
\begin{equation}
S_L^{-1}(s,T):=-(T^2-2{\rm Re}(s) T+|s|^2\mathcal{I})^{-1}(T-\overline{s}\mathcal{I}),\ \ \ s \in \rho_S(T)
\end{equation}
and
\begin{equation}
S_R^{-1}(s,T):=-(T-\overline{s}\mathcal{I})(T^2-2{\rm Re}(s) T+|s|^2\mathcal{I})^{-1},\ \ \ s \in  \rho_S(T),
\end{equation}
respectively.
Using the notion of $S$-spectrum and the notion of slice hyperholomorphic functions, see Section 4, we can define the quaternionic functional calculus, see \cite{formulation, JGA, book_functional}.
We point out that the $S$-resolvent operators are also used in Schur analysis in the realization of Schur functions in the slice hyperholomorphic setting see \cite{acs1, acs2, acs3}
and \cite{MR2002b:47144, adrs} for the classical case.

To set the framework in which we will work, we give some preliminaries.
Consider the complex plane $\mathbb{C}_I:=\mathbb{R}+I\mathbb{R}$, for $I\in \mathbb{S}$, where $\mathbb{S}$ is the unit sphere of purely imaginary quaternions. Observe that $\mathbb{C}_I$ can be identified with a complex plane since  $I^2=-1$ for every $I\in \mathbb{S}$.
Let $\Omega\subset \hh$ be a suitable domain that contains the $S$-spectrum of $T$.
We define the quaternionic functional calculus
  for  left slice hyperholomorphic functions $f:\Omega \to \hh$    as
\begin{equation}\label{quatinteg311def}
f(T)={{1}\over{2\pi }} \int_{\partial (\Omega\cap \mathbb{C}_I)} S_L^{-1} (s,T)\  ds_I \ f(s),
\end{equation}
where $ds_I=- ds I$;
for right slice hyperholomorphic functions, we define \begin{equation}\label{quatinteg311rightdef}
f(T)={{1}\over{2\pi }} \int_{\partial (\Omega\cap \mathbb{C}_I)} \  f(s)\ ds_I \ S_R^{-1} (s,T).
\end{equation}
These definitions are well posed since the integrals depend neither on the open set $\Omega$ nor on the complex plane
$\mathbb{C}_I$.
Using a similar idea, we define the projection operators which will provide the link between the spectral theorem and the $S$-spectrum.
\\
Our proofs will make use of a quaternionic version of Herglotz's theorem proved in the recent paper \cite{acks1}. This theorem will be the starting point to prove the quaternionic spectral theorem for unitary operators, in analogy with the classical case.

We have proved that if $U$ is a unitary operator acting on a quaternionic Hilbert space $\mathcal{H}$,
 then,  for $x$, $y\in \mathcal{H}$, there exists a spectral measure $E$ defined on the Borel sets of $[0,2\pi)$ such that for every slice continuous function $f\in \mathcal{S}(\sigma_S(U))$, we have
$$
\langle f(U) x, y \rangle = \int_0^{2\pi} f(e^{I t}) \langle dE(t)x, y \rangle, \quad\quad x,y, \in \cH.
$$
Moreover, for $t$ belonging to the Borel sets of $[0,2\pi)$, the measures
$$
\nu_{x,y}(t) = \langle E(t)x, y \rangle, \quad\quad x,y \in \cH,
$$
are related to the $S$-spectrum of $U$ by the quaternionic Riesz projectors through the relation
$$
\mathcal{P}(\sigma^0_S(U))=E(t_1)-E(t_0),
$$
where $\sigma^0_S(U)$ is the spectral set in the unit circle in  $\mathbb{C}_I$ delimited by the angles $t_0$, $t_1$.
\\
The plan of the paper is the following. In Section 2, we introduce the quaternionic Riesz projectors. Section 3 contains the proof of the main result of the paper, namely the spectral theorem for the unitary operatrs. In Section 4, we discuss the relation with the $S$-spectrum.

\section{Quaternionic Riesz projectors}
\setcounter{equation}{0}

In the following
we denote by  $\mathcal{B}(\mathcal{V})$ the space of bounded quaternionic linear operators on the left or on the right, since the results of this section hold in both cases.

The classical Riesz projectors are a powerful tool in spectral analysis and the  study of such projectors
is based on the resolvent equation.
Recently, in the paper \cite{acgs}, it has been shown that there exists a $S$-resolvent equation but in the quaternionic setting it involves both the $S$-resolvent operators. Precisely we have:
\begin{Tm}[The $S$-resolvent equation]\label{RLRESOLVEQ}
Let $T\in\mathcal{B}(\mathcal{V})$ and let $s$ and $p\in \rho_S(T)$. Then we have
\begin{equation}\label{RLresolv}\small
S_R^{-1}(s,T)S_L^{-1}(p,T)=((S_R^{-1}(s,T)-S_L^{-1}(p,T))p
-\overline{s}(S_R^{-1}(s,T)-S_L^{-1}(p,T)))(p^2-2s_0p+|s|^2)^{-1},
\end{equation}
but also
\begin{equation}\label{RLresolvII}\small
S_R^{-1}(s,T)S_L^{-1}(p,T)=(s^2-2p_0s+|p|^2)^{-1}(s(S_R^{-1}(s,T)-S_L^{-1}(p,T))
-(S_R^{-1}(s,T)-S_L^{-1}(p,T))\overline{p} ).
\end{equation}
\end{Tm}

The quaternionic functional calculus is  defined on the class of slice hyperholomorphic functions
$f: \Omega \subseteq\mathbb{H} \to \mathbb{H}$. Such functions have a Cauchy formula,
that works on specific domains which are called axially symmetric slice domain.
On this Cauchy formula is based the quaternionic functional calculus.

If we consider an element $I$ in the unit sphere of purely imaginary  quaternions
$$
\mathbb{S}=\{q=ix_1+jx_2+kx_3\ {\rm such \ that}\
x_1^2+x_2^2+x_3^2=1\}
$$
 then $I^2=-1$, and for this reason the elements of $\mathbb{S}$ are also called
imaginary units. Note that $\mathbb{S}$ is a 2-dimensional sphere in $\mathbb{R}^4$.
Given a nonreal quaternion $p=x_0+{\rm Im} (p)=x_0+I |{\rm Im} (p)|$, $I={\rm Im} (p)/|{\rm Im} (p)|\in\mathbb{S}$, we can associate to it the 2-dimensional sphere defined by
$$
[p]=\{x_0+I|{\rm Im} (p)|\ :\ I\in\mathbb{S}\}.
$$
 For any fixed $I\in\mathbb S$, the set $\mathbb C_I=\{u+Iv\ :\ u,v\in\mathbb R\}$ can be identified with the complex plane $\mathbb C$.

\begin{Dn}[Axially symmetric slice domain]\label{axsymm}
{\rm
Let $\Omega$ be a domain in $\mathbb{H}$.
We say that $\Omega$ is a
\textnormal{slice domain} (s-domain for short) if $\Omega \cap \mathbb{R}$ is
non empty and if $\Omega\cap \mathbb{C}_I$ is a domain in $\mathbb{C}_I$ for all $I \in \mathbb{S}$.
We say that $\Omega$ is
\textnormal{axially symmetric} if, for all $q \in \Omega$, the
$2$-sphere $[q]$ is contained in $\Omega$.}
\end{Dn}

\begin{Dn}
{\rm
An axially symmetric set $\sigma \subseteq \sigma_S(T)$  which is both open and closed in $\sigma_S(T)$ in its relative topology, is called a S-spectral set. For sake of simplicity
we will call it a spectral set.
}
\end{Dn}
The definition of a S-spectral set is suggested by the symmetry properties of the $S$-spectrum.
In fact, if $p\in \sigma_S(T)$,
  then all of the elements of the $2$-sphere $[p]$ are contained in  $\sigma_S(T)$.

\begin{Dn}
{\rm Let $T$ be a quaternionic linear operator on a quaternionic Banach space $\mathcal{V}$.
Denote by $\Omega_\sigma$  an axially symmetric s-domain that contains the spectral set $\sigma$
but not any other points of the $S$-spectrum.  Suppose that
the Jordan curves $\partial (\Omega_\sigma\cap \mathbb{C}_I)$ belong to the $S$-resolvent set $\rho_S(T)$, for every $I\in \mathbb{S}$.
We define the family $\mathcal{P}(\sigma)$  of quaternionic operators,  depending on the spectral sets $\sigma$, as
$$
\mathcal{P}(\sigma)=\frac{1}{2\pi }\int_{\partial (\Omega_\sigma\cap \mathbb{C}_I)}S_L^{-1}(s,T)ds_I.
$$}
\end{Dn}
The operators $\mathcal{P}(\sigma)$ are called (quaternionic) Riesz projectors.
\begin{Rk}\label{Rmk:2.5}{\rm
The definition of $\mathcal{P}(\sigma)$ can be given using the right $S$-resolvent operator $S_R^{-1}(s,T)$, that is
$$
\mathcal{P}(\sigma)=\frac{1}{2\pi }\int_{\partial (\Omega_\sigma\cap \mathbb{C}_I)}ds_IS_R^{-1}(s,T).
$$
}
\end{Rk}

Using the left $S$-resolvent operator we define the Riesz projectors associated with the $S$-spectrum. In \cite[Theorem 3.19]{acgs} we proved that $\mathcal{P}(\sigma)$ is a projector and that it commutes with $T$.
\\
The following lemma will be useful in the sequel.
\begin{La}\label{Lemma321gf}
Let $B\in \mathcal{B}(\mathcal{V})$ and let $\Omega$ be an axially symmetric s-domain.
\\
If  $p\in \Omega$, then
\begin{equation}\label{efdh}
\frac{1}{2\pi}\int_{\partial(\Omega\cap\mathbb{C}_I)}ds_I
(\overline{s}B-Bp)(p^2-2s_0p+|s|^2)^{-1}=B.
\end{equation}
Moreover, if $s\in \Omega$, then
\begin{equation}\label{efdh2}
\frac{1}{2\pi}\int_{\partial(\Omega\cap\mathbb{C}_I)}
(\overline{s}B-Bp)(p^2-2s_0p+|s|^2)^{-1}dp_I=-B.
\end{equation}
\end{La}
\begin{proof}
It follows the same lines of the proof of Lemma 3.18 in \cite{acgs}.
\end{proof}

\begin{Tm}\label{TEOREMAE}
Let $T$ be a quaternionic linear operator.
Then the family of operators $\mathcal{P}(\sigma)$ has the following properties{\rm :}
\begin{itemize}
\item[(i)]
$
(\mathcal{P}(\sigma))^2=\mathcal{P}(\sigma){\rm ;}
$
\item[(ii)]
$T\mathcal{P}(\sigma)=\mathcal{P}(\sigma)T{\rm ;}$
\item[(iii)]
$
\mathcal{P}(\sigma_S(T))=\mathcal{I}{\rm ;}
$
\item[(iv)]
$
\mathcal{P}(\emptyset)=0{\rm ;}
$
\item[(v)]
$
\mathcal{P}(\sigma\cup\delta)=\mathcal{P}(\sigma)+\mathcal{P}(\delta){\rm ;} \ \ \  \sigma \cap\delta=\emptyset,
$
\item[(vi)]
$
\mathcal{P}(\sigma\cap\delta)=\mathcal{P}(\sigma)\mathcal{P}(\delta).
$
\end{itemize}
\end{Tm}
\begin{proof}
Properties (i) and (ii) are proved in  Theorem 3.19 in \cite{acgs}.
Property (iii)
follows from the quaternionic functional calculus since
$$
T^m=\frac{1}{2\pi }\int_{\partial (\Omega\cap \mathbb{C}_I)}S_L^{-1}(s,T)ds_I \; s^m, \ \ \ m\in \mathbb{N}_0
$$
for $\sigma_S(T) \subset \Omega $, which for $m=0$ gives
$$
\mathcal{I}=\frac{1}{2\pi }\int_{\partial (\Omega\cap \mathbb{C}_I)}S_L^{-1}(s,T)ds_I.
$$
Property (iv) is a consequence of the functional calculus as well.

Property (v) follows from
\[
\begin{split}
\mathcal{P}(\sigma\cup\delta)&=\frac{1}{2\pi }\int_{\partial (\Omega_{\sigma\cup\delta}\cap \mathbb{C}_I)}S_L^{-1}(s,T)ds_I
\\
&
=\frac{1}{2\pi }\int_{\partial (\Omega_\sigma\cap \mathbb{C}_I)}S_L^{-1}(s,T)ds_I
+\frac{1}{2\pi }\int_{\partial (\Omega_\delta\cap \mathbb{C}_I)}S_L^{-1}(s,T)ds_I
\\
&=\mathcal{P}(\sigma)+\mathcal{P}(\delta).
\end{split}
\]

To prove (vi), assume that $\sigma \cap\delta\not=\emptyset$
and consider
\[
\begin{split}
\mathcal{P}(\sigma)\mathcal{P}(\delta)
&=
\frac{1}{(2\pi)^2 }
\int_{\partial (\Omega_\sigma\cap\mathbb{C}_I)}ds_IS_R^{-1}(s,T)
\int_{\partial (\Omega_\delta\cap \mathbb{C}_I)}S_L^{-1}(p,T)dp_I
\\
&=\frac{1}{(2\pi)^2 }\int_{ \partial ( \Omega_\sigma \cap \mathbb{C}_I) }  ds_I \int_{ \partial  (\Omega_\delta \cap \mathbb{C}_I) }
[S_R^{-1}(s,T)-S_L^{-1}(p,T)]
 p(p^2-2s_0p+|s|^2)^{-1}dp_I
\\
&
-\frac{1}{(2\pi)^2 }\int_{ \partial  (\Omega_\sigma \cap \mathbb{C}_I) }  ds_I
\int_{ \partial ( \Omega_\delta \cap \mathbb{C}_I) }\overline{s}[S_R^{-1}(s,T)-S_L^{-1}(p,T)]
(p^2-2s_0p+|s|^2)^{-1}
 dp_I,
 \end{split}
\]
where we have used the $S$-resolvent equation (see Theorem \ref{RLRESOLVEQ}). We rewrite the above relation as
\[
\begin{split}
\mathcal{P}(\sigma)\mathcal{P}(\delta) &=-\frac{1}{(2\pi)^2 }\int_{ \partial  (\Omega_\sigma \cap \mathbb{C}_I) }  ds_I \int_{ \partial  (\Omega_\delta \cap \mathbb{C}_I) }
[\overline{s}S_R^{-1}(s,T)-S_R^{-1}(s,T)p](p^2-2s_0p+|s|^2)^{-1}dp_I
\\
&
+\frac{1}{(2\pi)^2 }\int_{ \partial  (\Omega_\sigma \cap \mathbb{C}_I) }  ds_I
\int_{ \partial ( \Omega_\delta \cap \mathbb{C}_I) }[\overline{s}S_L^{-1}(p,T)-S_L^{-1}(p,T)p](p^2-2s_0p+|s|^2)^{-1}
 dp_I
 \\
& := \mathcal{J}_1+\mathcal{J}_2.
 \end{split}
\]
Now thanks to  Lemma \ref{Lemma321gf} and Remark \ref{Rmk:2.5}
we have
\[
\begin{split}
\mathcal{J}_1&=-\frac{1}{(2\pi)^2 }\int_{ \partial  (\Omega_\sigma \cap \mathbb{C}_I) }  ds_I \int_{ \partial ( \Omega_\delta \cap \mathbb{C}_I) }
[\overline{s}S_R^{-1}(s,T)-S_R^{-1}(s,T)p](p^2-2s_0p+|s|^2)^{-1}dp_I
\\
&
=\frac{1}{2\pi }\int_{ \partial  (\Omega_\sigma \cap \mathbb{C}_I) }  ds_IS_R^{-1}(s,T),\ \ {\rm for}\ \ \ s\in \Omega_\delta \cap \mathbb{C}_I
\\
&
=\frac{1}{2\pi }\int_{ \partial ( \Omega_\sigma \cap \mathbb{C}_I) } S_L^{-1}(s,T) ds_I,\ \  {\rm for}\ \ \ s\in \Omega_\delta \cap \mathbb{C}_I
\end{split}
\]
while $\mathcal{J}_1=0$ when $s\not\in \Omega_\delta \cap \mathbb{C}_I$ since
$$
\int_{ \partial ( \Omega_\delta \cap \mathbb{C}_I) }
[\overline{s}S_R^{-1}(s,T)-S_R^{-1}(s,T)p](p^2-2s_0p+|s|^2)^{-1}dp_I=0.
$$
Similarly, one can show that
\[
\begin{split}
\mathcal{J}_2
=\frac{1}{2\pi }
\int_{ \partial ( \Omega_\delta \cap \mathbb{C}_I) }S_L^{-1}(p,T)
 dp_I,\ \ \ \ {\rm for}\ \ \  p\in \Omega_\sigma \cap \mathbb{C}_I
\end{split}
\]
while $\mathcal{J}_2=0$ when $p\not\in \Omega_\sigma \cap \mathbb{C}_I$.
The integrals $\mathcal J_1$, $\mathcal J_2$ are either both zero or both nonzero,  so with a change of variable we get
$$
\mathcal{J}_1+\mathcal{J}_2=\frac{1}{2\pi }
\int_{ \partial ( \Omega_{\sigma\cap\delta} \cap \mathbb{C}_I) }S_L^{-1}(r ,T) dr _I=\mathcal{P}(\sigma\cap\delta).
$$
\end{proof}

From now on we will always work in quaternionic Hilbert spaces, so we will recall some definitions.

Let $\cH$ be a right linear quaternionic Hilbert space with an $\HH$-valued inner product $\langle \cdot, \cdot \rangle$ which satisfies
$$\langle x, y \rangle = \overline{ \langle y , x \rangle}{\rm ;}$$
$$\langle x, x \rangle \geq 0 \quad {\rm and} \quad \| x \|^2 := \langle x,x \rangle = 0 \Longleftrightarrow x = 0{\rm ;}$$
$$\langle x \alpha + y \beta, z \rangle = \langle x, z \rangle \alpha + \langle y, z \rangle \beta{\rm ;}$$
$$\langle x, y \alpha + z \beta \rangle = \overline{\alpha} \langle x, z \rangle + \overline{\beta} \langle x, z \rangle,$$
for all $x,y,z \in \cH$ and $\alpha,\beta \in \HH$. Any right linear quaternionic Hilbert space can be made also a left linear space, by fixing an Hilbert basis, see \cite{GMP}, Section 3.1.
We call an operator $A$ from the right quaternionic Hilbert space $\cH_1$, with inner product $\langle \cdot, \cdot \rangle_1$, to another right quaternionic Hilbert space $\cH_2$, with inner product $\langle \cdot, \cdot \rangle_2$, {\it right linear} if
$$A(x\alpha + y \beta) = (Ax)\alpha + (Ay)\beta,$$
for all $x,y$ in the domain of $A$ and $\alpha,\beta \in \HH$. We call an operator $A$ {\it bounded} if
$$\| A \| := \sup_{ \| x \| \leq 1 } \| A x \| < \infty.$$
Corresponding to any bounded right linear operator $A: \cH_1 \to \cH_2$ there exists a unique bounded right linear operator $A^*: \cH_2 \to \cH_1$ such that
$$\langle A x, y \rangle_2 = \langle x, A^* y \rangle_1,$$
and $\| A \| = \| A^* \|$ (see Proposition 6.2 in \cite{AlpayShapiro}).

Let $\cH$ be a right quaternionic Hilbert space with inner product $\langle \cdot, \cdot \rangle$. We call a right linear operator $U: \cH \to \cH$ {\it unitary} if
$$\langle U^* U x, y \rangle = \langle x,y \rangle,\quad\quad {\rm for}\; {\rm all} \; x,y \in \cH,$$
or, equivalently, $U^{-1} = U^*$.

\begin{Tm} Let $\cH$ be a right linear quaternionic Hilbert space and
let $U$ be a unitary operator on $\mathcal H$. Then the $S$-spectrum of $U$ belongs to the unit sphere of the quaternions.
\end{Tm}
\begin{proof}
See Theorem 4.8 in \cite{GMP}.
\end{proof}

By $ \cB([0,2\pi))$ we denote the Borel sets in  $[0,2\pi)$.
\begin{La} Let $x,y \in \mathcal{H}$ and let $\mathcal{P}(\sigma)$ be the projector associated with
 the unitary operator $U$. We define
$$
m_{x,y}(\sigma):= \langle \mathcal{P}(\sigma)x,y \rangle,\ \ \ \ x, \, y\in \cH,\ \  \ \sigma \in \cB([0,2\pi)).
$$
Then the $\HH$-valued measures $m_{x,y}$ defined on $\cB([0, 2\pi))$ enjoy the following properties
\begin{itemize}
\item[(i)] $m_{x \alpha +y \beta ,z} = m_{x,z}\alpha + m_{y,z} \beta${\rm ;}
\item[(ii)] $m_{x, y\alpha + z \beta} = \overline{\alpha} m_{x,y} + \overline{\beta} m_{x,z}${\rm ;}
\item[(iii)] $m_{x,y}([0, 2\pi)) \leq \| x \| \| y \|$,
\end{itemize}
where $x,y,z \in \cH$ and $\alpha,\beta \in \HH$.
\end{La}
\begin{proof} Properties (i) and (ii) follow from the properties of the quaternionic scalar product, while (iii) follows from Property (iii) in Theorem \ref{TEOREMAE} and the Cauchy-Schwarz inequality (see Lemma 5.6 in \cite{AlpayShapiro}).
\end{proof}

\section{The  spectral theorem for quaternionic unitary operators}
\setcounter{equation}{0}

We recall some classical results and also their quaternionic analogs which will be useful to prove a spectral theorem for quaternionic unitary operators.

\begin{Tm}[Herglotz's theorem]
\label{thm:Oct27yt1}
The function $n \mapsto r(n)$ from $\ZZ$ into $\CC^{s \times s}$ is positive definite if and only if there exists a unique $\CC^{s \times s}$-valued measure $\mu$ on $[0, 2\pi)$ such that
\begin{equation}
\label{eq:Oct27nv1}
r(n) = \int_0^{2\pi} e^{i n t} d\mu(t), \quad n \in \ZZ.
\end{equation}
\end{Tm}

Given $P \in \HH^{s \times s}$, there exist unique $P_1, P_2 \in \CC^{s \times s}$ such that $P = P_1 + P_2 j$. Recall the bijective homomorphism $\chi : \HH^{s \times s} \to \CC^{2s \times 2s}$ given by
\begin{equation}
\label{eq:Oct27jkl1}
\chi \hspace{0.5mm} P =  \begin{pmatrix} P_1 & P_2 \\ - \overline{P}_2 & \overline{P}_1 \end{pmatrix} \quad {\rm where} \; P = P_1 + P_2 j ,
\end{equation}

\begin{Dn}
\label{def:Oct27j1}
Given a $\HH^{s \times s}$-valued measure $\nu$, we may always write $\nu = \nu_1 + \nu_2 j$, where $\nu_1$ and $\nu_2$ are uniquely determined $\CC^{s \times s}$-valued measures. We call a measure $d\nu$ on $[0, 2\pi)$ {\it $q$-positive} if the $\CC^{2s \times 2s}$-valued measure
\begin{equation}
\label{eq:Nov7kr1}
\mu = \begin{pmatrix} \nu_1 & \nu_2 \\ \nu^*_2 & \nu_3 \end{pmatrix}, \quad {\rm where}\; \nu_3(t) = \nu_1(2\pi - t),\;\; t \in [0, 2\pi)
\end{equation}
is positive and, in addition,
$$\nu_2(t) = -\nu_2(2\pi -t)^T, \quad t \in [0, 2\pi).$$
\end{Dn}

\begin{Rk}
\label{rem:Oct27k1}{\rm
If $\nu$ is $q$-positive, then $\nu = \nu_1 + \nu_2 j$, where $\nu_1$ is a uniquely determined positive $\CC^{s \times s}$-valued measure and $\nu_2$ is a uniquely determined $\CC^{s \times s}$-valued measure.
}
\end{Rk}

\begin{Rk}
\label{rem:Nov7uuu1}
{\rm
If $r = (r(n))_{n \in \ZZ}$ is a $\HH^{s \times s}$-valued sequence on $\ZZ$ such that
$$r(n) = \int_0^{2\pi} e^{i n t} d\nu(t),$$
where $d\nu$ is a $q$-positive measure, then $r$ is Hermitian, i.e., $r(-n)^* = r(n)$.
}
\end{Rk}

The following result has been proved in \cite[Theorem 5.5]{acks1}.
\begin{Tm}[Herglotz's theorem for the quaternions]
\label{thm:Oct27j1}
The function $n \mapsto r(n)$ from $\ZZ$ into $\HH^{s \times s}$ is positive definite if and only if there exists a unique $q$-positive measure $\nu$ on $[0, 2\pi)$ such that
\begin{equation}
\label{eq:Oct27j1}
r(n) = \int_0^{2\pi} e^{i n t} d\nu(t), \quad n \in \mathbb{Z}.
\end{equation}
\end{Tm}

\begin{Rk}
\label{rem:Feb26yb1}
{\rm
For every $I \in \mathbb{S}$, there exists $J \in \mathbb{S}$ so that $I J = -J I$. Thus,
$\H = \mathbb C_I \oplus \mathbb C_I J$ and we may rewrite \eqref{eq:Oct27j1} as
\begin{equation}
\label{eq:Feb26um1}
r(n) = \int_0^{2\pi} e^{I n t} d\nu(t), \quad n \in \mathbb{Z},
\end{equation}
where $\nu = \nu_1 + \nu_2 J$ is a $q$-positive measure (in the sense that
$$\mu = \begin{pmatrix} \nu_1 & \nu_2 \\ \nu_2^* & \nu_3 \end{pmatrix}$$
is positive). Here $\nu_3(t) = \nu_1(2\pi - t)$.  }
\end{Rk}

For our purpose the scalar case will be important.

\begin{La}
\label{lem:Nov11a1}
If $U$ is a unitary operator on $\cH$, then $r_x = (r_x(n))_{n \in \ZZ}$, where $r_x(n) = \langle U^n x, x \rangle$ for $x \in \cH$, is an $\HH$-valued positive definite sequence.
\end{La}

\begin{proof}
If $\{ p_0, \ldots, p_N \} \subset \HH$, then
\begin{align*}
\sum_{m,n=0}^N \bar{p}_m r_x(n - m) p_n =& \; \sum_{m,n=0}^N \bar{p}_m \langle U^{n-m} x, x \rangle p_n \\
=& \; \sum_{m,n=0}^N \langle U^{n-m} x p_n, x p_m \rangle \\
=& \; \sum_{m,n=0}^N \langle U^n x p_n, U^m x p_m \rangle \\
=& \;  \langle \sum_{n=0}^N U^n x p_n, \sum_{m=0}^N U^m x p_m \rangle \\
=& \; \left\| \sum_{n=0}^N U^n x p_n \right\|^2 \geq 0.
\end{align*}
Thus, $r_x$ is a positive definite $\HH$-valued sequence.
\end{proof}

Let $r_x$ be as in Lemma \ref{lem:Nov11a1}. It follows from Theorem \ref{thm:Oct27j1} that there exists a unique $q$-positive measure $d\nu_x$ such that
\begin{equation}
\label{eq:Nov11kk1}
r_x(n) = \langle U^n x, x \rangle =  \int_0^{2\pi} e^{i n t } d\nu_x(t), \quad \quad n \in \ZZ.
\end{equation}
One can check that
\begin{align}
4 \langle U^n x, y \rangle =& \; \langle U^n (x+y), x+y \rangle - \langle U^n (x-y), x-y \rangle + i \langle U^n(x+yi), x+yi \rangle  \nonumber \\
& \; \;  - i \langle U^n (x - y i), x - yi \rangle + i \langle U^n (x - y j), x - y j \rangle k - i \langle U^n (x + y j) , x+ y j \rangle k \nonumber \\
& \; \; + \langle U^n(x+yk), x+yk \rangle k - \langle U^n(x-yk), x-yk \rangle k \label{eq:Nov11ki1}
\end{align}
and hence if we let
\begin{align}
4  \nu_{x,y} :=& \;  \nu_{x+y} -  \nu_{x-y} + i  \nu_{x+yi} - i  \nu_{x-yi} + i  \nu_{x-yj} k - i  \nu_{x+yj} k \nonumber \\
& \; \; +  \nu_{x+yk} k -  \nu_{x-yk} k, \label{eq:Nov11gkm1}
\end{align}
then
\begin{equation}
\label{eq:Nov11yj1}
\langle U^n x, y \rangle = \int_0^{2\pi} e^{i n t } d\nu_{x,y}(t), \quad\quad x,y \in \cH \quad {\rm and }\quad  n \in \ZZ.
\end{equation}

\begin{Tm}
\label{thm:Nov11kk1}
The $\HH$-valued measures $\nu_{x,y}$ defined on $\cB([0, 2\pi))$ enjoy the following properties{\rm :}
\begin{itemize}
\item[(i)] $\nu_{x \alpha +y \beta ,z} = \nu_{x,z}\alpha + \nu_{y,z} \beta,\ \ \alpha, \beta \in \mathbb{H}${\rm ;}
\item[(ii)] $\nu_{x, y\alpha + z \beta} = \bar{\alpha} \nu_{x,y} +\bar{\beta} \nu_{x,z},\ \ \alpha, \beta \in \mathbb{C}_i${\rm ;}
\item[(iii)] $\nu_{x,y}([0, 2\pi)) \leq \| x \| \| y \|$,
\end{itemize}
where $x,y,z \in \cH$ and $\alpha,\beta \in \HH$.
\end{Tm}

\begin{proof}
It follows from \eqref{eq:Nov11yj1} that
\begin{align*}
\int_0^{2\pi} e^{i n t} d\nu_{x\alpha + y \beta, z}(t) =& \; \langle U^n x, z \rangle \alpha + \langle U^n y, z \rangle \beta \\
=& \; \int_0^{2\pi} e^{i n t} (d \nu_{x,z}(t)\alpha  + d\nu_{y,z}(t) \beta), \quad\quad n \in \ZZ.
\end{align*}
Use the uniqueness in Theorem \ref{thm:Oct27j1} to conclude that
$$
\nu_{x\alpha + y \beta, z}(t) =  \nu_{x,z}(t)\alpha  + \nu_{y,z}(t) \beta
$$
 and hence we have proved (i). Property (ii) is proved in a similar fashion, observing that $\bar\alpha$, $\bar\beta$ commute with $e^{int}$.

If $n = 0$ in \eqref{eq:Nov11yj1}, then
$$\langle x , y \rangle = \int_0^{2\pi} d\nu_{x,y}(t) = \nu_{x,y}([0, 2\pi))$$
and thus we can use an analog of the Cauchy-Schwarz inequality (see Lemma 5.6 in \cite{AlpayShapiro}) to obtain
$$
\nu_{x,y}([0, 2\pi)) \leq \| x \| \| y \|$$
and hence we have proved (iii).
\end{proof}

\begin{Rk}
\label{rem:Nov11yy1}
{\rm
Contrary to the classical complex Hilbert space setting, $\nu_{x,y}$ need not equal ${\bar{\nu}_{y,x}}$.}
\end{Rk}

It follows from statements (i), (ii) and (iii) in Theorem \ref{thm:Nov11kk1} that $\phi(x) = \nu_{x,y}(\sigma)$, where $y \in \cH$ and $\sigma \in \cB([0, 2\pi))$ are fixed, is a continuous right linear functional. It follows from an analog of the Riesz representation theorem (see Theorem 6.1 in \cite{AlpayShapiro} or Theorem 7.6 in \cite{BDS}) that corresponding to any $x \in \cH$, there exists a uniquely determined vector $w \in \cH$ such that
$$\phi(x) = \langle x, w \rangle,$$
i.e. $\nu_{x,y}(\sigma) = \langle x, w \rangle$. Use (i) and (ii) in Theorem \ref{thm:Nov11kk1} to deduce that $w = E(\sigma)^* y$. The uniqueness of $E$ follows readily from the construction. Thus, we have
\begin{equation}
\label{eq:Nov11uujj1}
\nu_{x,y}(\sigma) = \langle E(\sigma)x, y \rangle, \quad\quad x,y \in \cH \quad {\rm and} \quad \sigma \in \cB([0,2\pi)),
\end{equation}
whence
\begin{equation}
\label{eq:Nov11jeje1}
\langle U^n x, y \rangle = \int_0^{2\pi} e^{ i n t} \langle dE(t) x, y \rangle.
\end{equation}

To prove the main properties of the operator $E$ we need a uniqueness results on quaternionic measures
which is a corollary of the following:
\begin{Tm}
\label{thm:Feb23yn1}
Let $\mu$ and $\nu$ be $\CC$-valued measures on $[0, 2\pi)$. If
\begin{equation}
\label{eq:Feb23udn1}
r(n) = \int_0^{2\pi} e^{i n t} d\mu(t) = \int_0^{2\pi} e^{i n t} d\nu(t), \quad n \in \mathbb{Z},
\end{equation}
then $\mu = \nu$.
\end{Tm}
\begin{proof}
See, e.g., Theorem 1.9.5 in \cite{Sasvari}.
\end{proof}

\begin{Tm}
\label{thm:Feb23u1}
Let $\mu$ and $\nu$ be $\HH$-valued measures on $[0, 2\pi)$. If
\begin{equation}
\label{eq:Feb23un1}
r(n) = \int_0^{2\pi} e^{i n t} d\mu(t) = \int_0^{2\pi} e^{i n t} d\nu(t), \quad n \in \ZZ,
\end{equation}
then $\mu = \nu$.
\end{Tm}

\begin{proof}
Write $r(n) = r_1(n) + r_2(n)j$, $\mu = \mu_1 + \mu_2j$ and $\nu = \nu_1 + \nu_2 j$, where $r_1(n), r_2(n) \in \CC$ and $\mu_1, \mu_2, \nu_1, \nu_2$ are $\CC$-valued measures on $[0, 2\pi)$. It follows from \eqref{eq:Feb23un1} that
$$r_1(n) = \int_0^{2\pi} e^{i n t} d\mu_1(t) = \int_0^{2\pi} e^{i n t} d\nu_1(t), \quad n \in \ZZ$$
and
$$r_2(n) = \int_0^{2\pi} e^{i n t} d\mu_2(t) = \int_0^{2\pi} e^{i n t} d\nu_2(t), \quad n \in \ZZ.$$
Use Theorem \ref{thm:Feb23yn1} to conclude that $\mu_1 = \nu_1$, $\mu_2 = \nu_2$ and hence that $\mu = \nu$.
\end{proof}

\begin{Tm}
\label{thm:Nov11ikm1}
The operator $E$ given in \eqref{eq:Nov11uujj1} enjoys the following properties{\rm :}
\begin{itemize}
\item[(i)] $\| E(\sigma) \| \leq 1${\rm ;}
\item[(ii)] $E(\emptyset) = 0$ and $E([0, 2\pi)) = I_{\cH}${\rm ;}
\item[(iii)] If $\sigma \cap \tau = \emptyset$, then $E( \sigma \cup \tau ) = E(\sigma) + E(\tau)${\rm ;}
\item[(iv)] $E(\sigma \cap \tau) = E(\sigma) E(\tau)${\rm ;}
\item[(v)] $E(\sigma)^2 = E(\sigma)${\rm ;}
\item[(vi)] $E(\sigma)$ commutes with $U$ for all $\sigma \in \cB([0, 2\pi))${\rm .}
\end{itemize}
\end{Tm}

\begin{proof}
Use \eqref{eq:Nov11uujj1} with $y = E(\sigma) x$ and (iii) in Theorem \eqref{thm:Nov11kk1} to obtain
$$\| E(\sigma) x \|^2 \leq \| x \| \| E(\sigma) x \|,$$
whence we have shown (i). The first part of property (ii) follows directly from the fact that $\nu_{x,y}(\emptyset) = 0$. The last part follows from \eqref{eq:Nov11jeje1} when $n = 0$.
Statement (iii) follows easily from the additivity of the measure $\nu_{x,y}$.

We will now prove property (iv). It follows from \eqref{eq:Nov11jeje1} that
\begin{align*}
\langle U^{n+m} x, y \rangle =& \; \int_0^{2\pi} e^{i n t} e^{i m t} \langle dE(t)x, y \rangle \\
=& \; \langle U^n (U^m x), y \rangle \\
=& \; \int_0^{2\pi} e^{i n t} d\langle E(t) U^m x, y \rangle.
\end{align*}
Using the uniqueness in Theorem \ref{thm:Feb23u1} we obtain
$$e^{i m t} d\langle E(t) x, y \rangle = \langle dE(t) U^m x, y \rangle$$
and hence, denoting by $\mathbf{1}_\sigma$ the characteristic function of the set $\sigma$, we have
$$\int_0^{2\pi} \mathbf{1}_\sigma(t) e^{i m t} \langle dE(t) x, y \rangle = \langle E(\sigma) U^m x , y \rangle.$$
 But
$$\int_0^{2\pi} \mathbf{1}_\sigma(t) e^{i m t} \langle dE(t) x, y \rangle = \langle U^k x, E(\sigma)^* y \rangle = \int_0^{2\pi} e^{i m t} d\langle E(t) x, E(\sigma)^* y \rangle.$$
Using the uniqueness in Theorem \ref{thm:Feb23u1} once more we get
$$\mathbf{1}_\sigma(t) d\langle E(t)x, y \rangle =  \langle dE(t) x, E(\sigma)^* y \rangle$$
and hence
$$\int_0^{2\pi} \mathbf{1}_\tau(t) \mathbf{1}_\sigma(t)  \langle dE(t) x, y \rangle = \langle E(t)x, E(\sigma)^* y \rangle$$
and thus
$$ \langle E( \sigma \cap \tau ) x, y \rangle = \langle E(\sigma) E(\tau) x, y \rangle.$$
Property (v) is obtained from (iv) by letting $\sigma = \tau$.

Finally, since $U$ is unitary one can check that
$$\langle U(x \pm U^* y), x \pm U^* y \rangle = \langle U(Ux \pm y), Ux \pm y \rangle$$
and hence from \eqref{eq:Nov11yj1} and the uniqueness in Theorem \ref{thm:Feb23u1} we obtain $\nu_{x \pm U^*y} =  \nu_{Ux \pm y}$.
Similarly,
$$\nu_{x \pm U^* y i} = \nu_{Ux \pm yi}$$
$$\nu_{x\pm U^* y j} = \nu_{U x \pm yj }$$
and
$$\nu_{x \pm U^* y k} = \nu_{Ux \pm yk}.$$
It follows from \eqref{eq:Nov11gkm1} that
$$\nu_{x, U^* y} = \nu_{Ux, y}.$$
Now use \eqref{eq:Nov11uujj1} to obtain
$$\langle E(\sigma) x, U^* y \rangle = \langle E(\sigma) U x, y \rangle,$$
i.e.,
$$ \langle U E(\sigma) x, y \rangle = \langle E(\sigma) U x , y \rangle,\quad\quad x,y \in \cH.$$
\end{proof}

Given any quaternionic Hilbert space $\cH$, there exists a subspace $\mathcal{M} \subset \cH$ on $\CC$ so that for any $x \in \cH$ we have
$$x = x_1 + x_2 j, \quad x_1,x_2 \in \mathcal{M}.$$

\begin{Tm}
Let $U$ be a unitary operator on a quaternionic Hilbert space $\cH$ and let $E$ be the corresponding operator given by \eqref{eq:Nov11uujj1}. $E$ is self-adjoint if and only if $U: \mathcal{M} \to \mathcal{M}$, where $\mathcal{M}$ is as above.
\end{Tm}

\begin{proof}
If $E = E^*$, then it follows from \eqref{eq:Nov11uujj1} that $\nu_{x,y} = \bar{\nu}_{y,x}$ for all $x,y \in \cH$. In particular, we get $\nu_{x,x} = \bar{\nu}_{x,x}$, i.e.
\begin{equation}
\label{eq:Feb18au1}
\nu_x = \bar{\nu}_x, \quad x \in \cH.
\end{equation}
Since $\nu_x$ is a $q$-positive measure we may write $\nu_x = \alpha_x + \beta_x j$, where $\alpha_x$ is a positive Borel measure on $[0, 2\pi)$ and $\beta_x$ is a complex Borel measure on $[0, 2\pi)$. It follows from \eqref{eq:Feb18au1} that
$$\beta_x = - \beta_x,$$
i.e. $\beta_x = 0$. Thus, we may make use of the spectral theorem for unitary operators on a complex Hilbert space (see, e.g., Section 31.7 in \cite{Lax}) to deduce that $U: \mathcal{M} \to \mathcal{M}$. Conversely, if $U: \mathcal{M} \to \mathcal{M}$, then the spectral theorem for unitary operators on a complex Hilbert space yields that $E = E^*$.
\end{proof}

If $U: \HH^n \to \HH^n$ is unitary, then \eqref{eq:Nov11jeje1} and Theorem \ref{thm:Nov11ikm1} assert that
\begin{equation}
\label{eq:Feb20yn1m}
U = \sum_{a=1}^n e^{i \theta_a} P_a,
\end{equation}
where $\theta_1, \ldots, \theta_n \in [0, 2\pi)$ and $P_1, \ldots, P_n$ are oblique projections (i.e. $(P_a)^2 = P_a$ but $(P_a)^*$ need not equal $P_a$). Corollary 6.2 in Zhang \cite{Zhang} asserts, in particular, the existence of $V: \HH^n \to \HH^n$ which is unitary and $\theta_1, \ldots, \theta_n \in [0, 2\pi)$ so that
\begin{equation}
\label{eq:Feb18unb1}
U = V^* {\rm diag}(e^{i \theta_1}, \ldots, e^{i \theta_n}) V.
\end{equation}
In the following remark we will explain how \eqref{eq:Feb20yn1m} and \eqref{eq:Feb18unb1} are consistent.

\begin{Rk}
\label{rem:Feb18z1}
{\rm Let $U: \HH^n \to \HH^n$ be unitary. Let $V$ and $\theta_1, \ldots, \theta_n$ be as above. If we let $e_a = (0, \ldots, 0, 1, 0, \ldots, 0)^T \in \HH^n$, where the $1$ is the $a$-th position, then we can rewrite \eqref{eq:Feb18unb1} as
$$
U = \sum_{a=1}^n V^* e^{i \theta_a} e_a e_a^* V.
$$
Note that $V^* e^{i \theta_a} e_a e_a^* V = e^{i \theta_a} V^* e_a e_a^* V$ if and only if $V: \CC^n \to \CC^n$. In this case $U: \CC^n \to \CC^n$ and
$$U = \sum_{a=1}^n e^{i \theta_a} P_a,$$
where $P_a$ denotes the orthogonal projection given by $V^* e^{i \theta_a} e_a e_a^* V$.
}
\end{Rk}

\begin{Rk}
\label{SpectrIJK}{\rm
Observe that in the proof of the spectral theorem for $U^n$ we have taken the imaginary units $i$, $j$, $k$ for  the quaternions and we have determined spectral measures $\langle dE(t) x, y \rangle$ that are supported on the unit circle in $\mathbb{C}_i$.
In the case one uses other orthogonal units $I$, $J$ and $K\in \mathbb{S}$ to represent quaternions,
then the spectral measures are supported on the unit circle  in $\mathbb{C}_I$.
}
\end{Rk}

 Observe that \eqref{eq:Nov11jeje1} provides a vehicle to define a functional calculus for unitary operators on a quaternionic Hilbert space. For a continuous $\HH$-valued function $f$ on the unit circle, which will be approximated by the polynomials $\sum_{k}e^{i k t} a_k $.
We will consider a subclass of continuous quaternionic-valued functions  defined as follows, see \cite{GMP}:
\begin{Dn}\label{defdiN_c}
{\rm
The quaternionic linear space of slice continuous functions on an axially symmetric subset $\Omega$ of $\mathbb H$, denoted by $\mathcal{S}(\Omega)$ consists of functions of the form $f(u+Iv)=\alpha(u,v)+I\beta(u,v)$ where $\alpha,\beta$ are quaternionic valued functions such that $\alpha(x,y)=\alpha(u,-v)$, $\beta(u,v)=-\beta(u,-v)$ and $\alpha$, $\beta$ are continuous functions. When $\alpha,\beta$ are real valued we say that the continuous slice function is intrinsic. The subspace of intrinsic continuous slice functions is denoted by $\mathcal{S}_{\mathbb R}(\Omega)$.}
\end{Dn}
It is important to note that any polynomial of the form $P(u+Iv)=\sum_{k=0}^n (u+Iv)^n a_n$, $a_n\in\mathbb H$ is a slice continuous function in the whole $\mathbb H$. A trigonometric polynomial of the form $P(e^{It})=\sum_{m=-n}^n e^{Imt} a_m$ is a slice continuous function on $\partial \mathbb B$, where $\mathbb B$ denotes the unit ball of quaternions.\\
\\
Let us now denote by $\mathcal{PS}(\sigma_S(T))$ the set of slice continuous functions $f(u+Iv)=\alpha(u,v)+I\beta(u,v)$ where $\alpha$, $\beta$ are polynomials in the variables $u,v$.
\\
In the sequel we will work on the complex plane $\mathbb C_I$ and we denote by $\mathbb{T}_I$ the boundary of $\mathbb B\cap\mathbb C_I$.
Any other choice of an imaginary unit in the unit sphere $\mathbb S$ will provide an analogous result.

\begin{Rk}
\label{rem:Feb27qa1}
{\rm For every $I \in \mathbb{S}$, there exists $J \in \mathbb{S}$ so that $I J = -J I$. Bearing in mind Remark \ref{rem:Feb26yb1}, we can construct $\nu_{x,y}^{(J)}$ so that \eqref{eq:Nov11yj1} can also be written as
\begin{equation}
\label{eq:Feb26ub1}
\langle U^n x, y \rangle = \int_0^{2\pi} e^{I n t } d\nu^{ (J)}_{x,y}(t), \quad\quad x,y \in \cH \quad {\rm and }\quad  n \in \ZZ.
\end{equation}
Consequently, \eqref{eq:Nov11jeje1} can be written as
\begin{equation}
\label{eq:Feb27unq1}
\langle U^n x, y \rangle = \int_0^{2\pi} e^{i n t} \langle E_J(t) x, y \rangle,
\end{equation}
where $E_J$ is given by
$$\nu_{x,y}^{(J)}(\sigma) = \langle E_J(\sigma)x, y \rangle, \quad x,y \in \mathcal{H} \quad {\rm and} \quad \sigma \in {\rm B}(\mathbb{T}_I).$$
Moreover, $E_J$ satisfy properties (i)-(v) listed in Theorem \ref{thm:Nov11ikm1}.
}

\end{Rk}

\begin{Tm}[The spectral theorem for quaternionic unitary operators]
\label{fctcalc}
Let $U$ be an unitary operator on a right linear quaternionic Hilbert space $\mathcal{H}$. Let $I, J \in \mathbb{S}$, $I$ orthogonal to $J$. Then there exists a unique spectral measure $E_J$ defined on the Borel sets of $\mathbb{T}_I$ such that for every slice continuous function $f\in \mathcal{S}(\sigma_S(U))$, we have
$$
f(U) = \int_0^{2\pi} f(e^{I t}) dE_J(t).
$$
\end{Tm}
\begin{proof}
Let us consider a polynomial $P(t) = \sum_{m=-n }^n  e^{I m t}a_m$ defined on $\mathbb{T}_I$. Then using \eqref{eq:Feb27unq1} we have
$$
\langle U^m x, y \rangle = \int_0^{2\pi} e^{Im t} \langle dE_J(t)x, y \rangle \quad\quad x,y, \in \cH.
$$
By linearity, we can define
$$
\langle P(U) x, y \rangle = \int_0^{2\pi} P(e^{I t}) \langle dE_J(t)x, y \rangle,\quad\quad x,y, \in \cH.
$$

The map $\Psi:\ \mathcal{PS}(\sigma_S(U))\to \mathcal H$ defined by $\psi_U(P) = P(U)$ is $\mathbb{R}$-linear.
 By fixing a basis for $\mathbb H$, e.g. the basis $1,i,j,k$, each slice continuous function $f$ can be decomposed using intrinsic function, i.e. $f=f_0+f_1i+f_2j+f_3k$ with $f_\ell\in\mathcal{S}_{\mathbb R}(\sigma_S(U))$, $\ell=0,\ldots,3$, see \cite[Lemma 6.11]{GMP}. For these functions the spectral mapping theorem holds, thus $f_\ell(\sigma_S(U))=\sigma_S(f_\ell (U))$ and so $\| f_\ell(U)\| =\| f_\ell \|_\infty$, see \cite[Theorem 7.4]{GMP}.
 The map $\psi$ is continuous and so there exists  $C>0$, that does not depend on $\ell$,  such that
$$
\| P(U) \|_{\mathcal{H}}  \leq   C \max_{ t \in \sigma_S(U) } |P(t) |.
$$
A slice continuous function $f\in\mathcal S(\sigma_S(U))$ is defined on an axially symmetric subset $K\subseteq \mathbb T$ and thus it can be written as a function
$f(e^{It})=\alpha(\cos t,\sin t) +I \beta (\cos t,\sin t)$. By fixing a basis of $\mathbb H$, e.g. $1,i,j,k$, $f$ can be decomposed into four slice continuous  intrinsic functions $f_\ell(\cos t,\sin t)=\alpha_\ell (\cos t,\sin t) +I \beta_\ell (\cos t,\sin t)$, $\ell=0,\ldots,3$, such that $f= f_0+f_1i+f_2j+f_3k$.

By the  Weierstrass approximation theorem for trigonometric polynomials, see, e.g., Theorem 8.15 in \cite{Rudin}, each function $f_\ell$
can be approximated by a sequence of polynomials $$\tilde R_{\ell n} = \tilde a_{\ell n} (\cos t,\sin t) +I \tilde b_{\ell n} (\cos t,\sin t),$$ $\ell=0,\ldots,3$ which
tend uniformly to $f_\ell$. These polynomials do not necessarily belong to the class of the continuous slice functions since $\tilde a_{\ell n}, \tilde b_{\ell n}$ do not satisfy, in general, the even and odd conditions in Definition \ref{defdiN_c}. However, by setting
$$a_{\ell n}(u,v)= \frac 12(\tilde a_{\ell n}(u,v)+\tilde a_{\ell n}(u,-v)),$$
$$b_{\ell n}(u,v)= \frac 12(\tilde b_{\ell n}(u,-v)-\tilde b_{\ell n}(u,v))$$
we obtain that the sequence of polynomials
$a_{\ell n}+Ib_{\ell n}$ still converges to $f_\ell$, $\ell=0,\ldots, 3$. By putting $R_{\ell n} = a_{\ell n} (\cos t,\sin t) +I b_{\ell n} (\cos t,\sin t)$, $\ell=0,\ldots,3$ and $R_n=R_{0n}+R_{1n}i+R_{2n}j+R_{3n}k$ we have a sequence of slice continuous polynomials $R_n(e^{It})$ converging to $f(e^{It})$ uniformly on $\mathbb R$.

By the previous discussion, $\{R_n(U)\}$ is a Cauchy sequence in the space of bounded linear operators since
$$
\| R_n(U)-R_m(U) \| \leq  C \max_{ t \in \sigma_S(U) } |R_n(t)-R_m(t) |,
$$
so  as ${R_n(U)}$  has a limit which we denote by $f(U)$.
\end{proof}

\begin{Rk}
\label{rem:Feb26un1}
{\rm
Fix $I \in \mathbb{S}$. It is worth pointing out that $f(u+Iv) = (u+Iv)^{-1}$ is an intrinsic function on $\C_{I} \cap \partial \mathbb{B}$, where
$\partial \mathbb{B} = \{ q \in \H: |q| = 1 \}$, since
$$f(u+Iv) = \frac{u}{u^2+v^2} +  \left(\frac{-v}{u^2+v^2}\right)J.$$
Thus, using Theorem \ref{fctcalc}, we may write
\begin{equation}
\label{eq:Feb27unb1}
U^{-1} = \int_0^{2\pi} e^{-It} d E_J(t)
\end{equation}
and
\begin{equation}
U = \int_0^{2\pi} e^{It} dE_J(t).
\end{equation}
}
\end{Rk}

\section{The $S$-spectrum and the spectral theorem }

We now want to show that the spectral theorem is based on the $S$-spectrum.
We will be in need of the  Cauchy formula for slice hyperholomorphic functions, see \cite{book_functional} for more details.

\begin{Dn}[Cauchy kernels]
{\rm We define the  (left) Cauchy kernel, for $q\not\in[s]$, by
\begin{equation}\label{esse-1Left}
S_L^{-1}(s,q):=-(q^2-2q{\rm Re}(s)+|s|^2)^{-1}(q-\bar s).
\end{equation}
We define the  right Cauchy kernel, for $q\not\in[s]$, by
\begin{equation}\label{esseR}
S^{-1}_R(s,q):=-(q-\bar s)(q^2-2{\rm Re}(s)q+|s|^2)^{-1}.
\end{equation}}
\end{Dn}
\begin{Tm}\label{Cauchynuovo}\index{Cauchy!integral formula}
Let $\Omega \subseteq \mathbb{H}$ be an axially symmetric s-domain  such that
$\partial (\Omega \cap \mathbb{C}_I)$ is union of a finite number of
continuously differentiable Jordan curves, for every $I\in\mathbb{S}$. Let $f$ be
a slice hyperholomorphic function on an open set containing $ \overline{\Omega}$ and, for any $I\in \mathbb{S}$,  set  $ds_I=-Ids$.
Then for every $q=u+I_qv\in \Omega$ we have{\rm :}
\begin{equation}\label{integral}
 f(q)=\frac{1}{2 \pi}\int_{\partial (\Omega \cap \mathbb{C}_I)} S_L(s,q) ds_I f(s).
\end{equation}
Moreover,
the value of the integral depends neither on $\Omega$ nor on the  imaginary unit
$I\in\mathbb{S}$.

\noindent
If $f$ is a right slice regular function on a set that contains $\overline{\Omega}$,
then
\begin{equation}\label{Cauchyright_qua}
 f(q)=\frac{1}{2 \pi}\int_{\partial (\Omega\cap \mathbb{C}_I)}  f(s)ds_I S_R^{-1}(s,q).
 \end{equation}
Moreover,
the value of the integral depends neither on $\Omega$ nor on the  imaginary unit
$I\in\mathbb{S}$.
\end{Tm}
We conclude the paper with the following result, based on the Cauchy formula, that shows the relation between the spectral measures and the $S$-spectrum.
\begin{Tm}\label{specmeasures}
 Fix $I, J \in \mathbb{S}$, with $I$ orthogonal to $J$.
Let $U$ be an unitary operator on a right linear quaternionic Hilbert space $\mathcal{H}$ and let $E(t)=E_J(t)$ be its spectral measure.
Assume that $\sigma^0_S(U)\cap\mathbb{C}_I$ is contained in  the arc of the unit circle in $\mathbb C_I$ with endpoints
$t_0$ and $t_1$.
Then
$$
\mathcal{P}(\sigma^0_S(U))=E(t_1)-E(t_0).
$$
\end{Tm}
\begin{proof}
The spectral theorem implies that the operator $S_R^{-1}(s,U)$ can be written as
$$
S_R^{-1}(s,U)=\int_0^{2\pi}S_R^{-1}(e^{It},s) dE(t).
$$
The Riesz projector is given by
$$
\mathcal{P}(\sigma^0_S(U))=\frac{1}{2\pi}\int_{\partial(\Omega_0\cap\mathbb{C}_I)} ds_I S_R^{-1}(s,U)
$$
where $\Omega_0$ is an open set containing $\sigma^0_S(U)$ and such that $\partial(\Omega_0\cap\mathbb{C}_I)$ is a smooth closed curve in $\mathbb{C}_I$.
Now we write
$$
\mathcal{P}(\sigma^0_S(U))=\frac{1}{2\pi}\int_{\partial(\Omega_0\cap\mathbb{C}_I)} ds_I \Big(\int_0^{2\pi}S_R^{-1}(e^{It},s) dE(t) \Big)
$$
and using the Fubini theorem we get
$$
\mathcal{P}(\sigma^0_S(U))=\int_0^{2\pi}  \Big( \frac{1}{2\pi}\int_{\partial(\Omega_0\cap\mathbb{C}_I)} ds_I S_R^{-1}(e^{It},s) \Big) dE(t).
$$
It follows from the Cauchy formula that
$$
\frac{1}{2\pi}\int_{\partial(\Omega_0\cap\mathbb{C}_I)} ds_I S_R^{-1}(e^{It},s)=\mathbf{1}_{[t_0,t_1]},
$$
where $\mathbf{1}_{[t_0,t_1]}$ is the characteristic function of the set $[t_0,t_1]$, and so we get the statement, since
$$
\mathcal{P}(\sigma^0_S(U))=\int_0^{2\pi}  \mathbf{1}_{[t_0,t_1]}dE(t)=E(t_1)-E(t_2).
$$

\end{proof}

\end{document}